\documentclass[11pt]{article}
\newcommand{\documentdate}{31 X 2024}

\usepackage{a4wide,latexsym,amsmath,varioref,graphicx,xcolor,amssymb,diagbox}

\title{Complexity of Adagrad and other first-order methods\\
  for nonconvex optimization problems \\
  with bounds constraints}

\author{
   S. Gratton%
   \thanks{Universit\'e de Toulouse, INP, IRIT, Toulouse, France. Email:
     serge.gratton@enseeiht.fr. Work partially supported by 3IA Artificial and
     Natural Intelligence Toulouse Institute (ANITI), French "Investing for the Future
     - PIA3" program under the Grant agreement ANR-19-PI3A-0004"}, 
   ~S. Jerad%
   \thanks{ANITI, Universit\'e de Toulouse, INP, IRIT, Toulouse, France. Email:
     sadok.jerad@enseeiht.fr}
   ~and Ph. L. Toint%
   \thanks{NAXYS, University of Namur, Namur, Belgium. Email:
     philippe.toint@unamur.be.
     Partially supported by ANITI.}
}

\newtheorem{theorem}{Theorem}[section]
\newtheorem{lemma}[theorem]{Lemma}
\newtheorem{corollary}[theorem]{Corollary}
\newcommand{\numsection}[1]{\section{#1}\setcounter{equation}{0}}
\newcommand{\appnumsection}[1]{\section*{#1}
  \renewcommand{\theequation}{A.\arabic{equation}}
  \renewcommand{\thetheorem}{A.\arabic{theorem}}
  \renewcommand{\thetable}{A.\arabic{table}}
  \renewcommand{\thefigure}{A.\arabic{figure}}
  \renewcommand{\thesection}{A} }
\renewcommand{\theequation}{\arabic{section}.\arabic{equation}}

\newcounter{algo}[section]
\renewcommand{\thealgo}{\thesection.\arabic{algo}}
\newcommand{\llem}[2]{\vspace{\baselineskip} 
\noindent\framebox[\textwidth]{\parbox{0.95\textwidth}{
\begin{lemma} \label{#1} \rm #2 \end{lemma} } } \vspace{\baselineskip} }

\newcommand{\algo}[3]{\refstepcounter{algo}
\begin{center}\begin{figure}[htbp]
\framebox[\textwidth]{
\parbox{0.95\textwidth} {\vspace{\topsep}
{\bf Algorithm \thealgo : #2}\label{#1}\\
\vspace*{-\topsep} \mbox{ }\\
{#3} \vspace{\topsep} }}
\end{figure}\end{center}}
\newcommand{\bpr}{{\bf Proof.} \hspace{1.5mm}}
\newcommand{\epr}{\hfill $\Box$ \vspace*{1em}}
\newcommand{\proof}[1]{
\begin{list}{}{
\setlength{\topsep}{0.0pt}
\setlength{\partopsep}{0.0pt}
\setlength{\leftmargin}{0.025\textwidth}
\setlength{\rightmargin}{0.5\leftmargin}
\setlength{\labelwidth}{0.5\leftmargin}
\setlength{\labelsep}{0.25\leftmargin}}
\item \bpr #1 \epr \noindent
\end{list}}
\newcommand{\lthm}[2]{\vspace{\baselineskip} 
\noindent\framebox[\textwidth]{\parbox{0.95\textwidth}{
\begin{theorem} \label{#1} \rm #2 \end{theorem} } } \vspace{\baselineskip} }
\newcommand{\beqn}[1]{\begin{equation}\label{#1}}
\newcommand{\eeqn}{\end{equation}}
\newcommand{\req}[1]{(\ref{#1})}
\newcommand{\ms}{\;\;\;\;}
\newcommand{\tim}[1]{\;\; \mbox{#1} \;\;}
\newcommand{\ii}[1]{\{ 1, \ldots, #1 \}}
\newcommand{\iiz}[1]{\{ 0, \ldots, #1 \}}
\newcommand{\iibe}[2]{\{ #1, \ldots, #2 \}}
 
\newcommand{\calF}{{\cal F}}

\newcommand{\calN}{{\cal N}} 
\newcommand{\calO}{{\cal O}}

\renewcommand{\Re}{\hbox{I\hskip -2pt R}}
\newcommand{\smallRe}{\hbox{\footnotesize I\hskip -2pt R}}
\newcommand{\bigfrac}[2]{\frac{\displaystyle #1}{\displaystyle #2}}
\newcommand{\bigsum}{\displaystyle \sum}

\newcommand{\sfrac}[2]{{\scriptstyle \frac{#1}{#2}}}
\newcommand{\kap}[1]{\kappa_{\mbox{\tiny #1}}}
\newcommand{\eqdef}{\stackrel{\rm def}{=}}
\newcommand{\al}[1]{{\footnotesize{\sf #1}}}

\newcommand{\tal}[1]{{\normalsize {\sf #1}}}
\newcommand{\half}{\sfrac{1}{2}}

\newcommand{\flow}{f_{\rm low}}
\newcommand{\sgn}{{\rm sign}}
\newcommand{\neol}{\nonumber\\}
\newcommand{\kB}{\kap{B}}
\DeclareMathOperator*{\average}{average}

\newcommand{\comment}[1]{}

\date{\documentdate}

\begin{document}

\maketitle

\begin{abstract}
A parametric class of trust-region algorithms for constrained
nonconvex optimization is analyzed, where the objective function is
never computed. By defining appropriate first-order stationarity
criteria, we are able to extend the Adagrad method to the newly
considered problem and retrieve the standard complexity rate of the
projected gradient method that uses both the gradient and objective
function values. Furthermore, we propose an additional
iteration-dependent scaling with slightly inferior theoretical
guarantees. In both cases, the bounds are essentially sharp, and
curvature information can be used to compute the stepsize. Initial
experimental results for noisy bound-constrained instances illustrate
the benefits of the objective-free approach.
\end{abstract} 


\numsection{Introduction}

It is a truism to say that adaptive-gradient methods such as Adagrad,
ADAM and cousins are at the heart of machine learning algorithms, and
the literature covering them is vast (see
\cite{LecuBottBengHaff98,DuchHazaSing11,TielHint12,KingBa15,
  WardWuBott19,ReddKaleKuma18,DefoBottBachUsun22} to cite only a few
significant contributions). In line with current technology for neural
network training, most of these minimization methods have been
considered for nonconvex \textit{unconstrained} problems.  While this
is adequate for a broad category of applications, it is however
difficult to use them in a context where a priori information on the
problem at hand is available, often in the form of
constraints. Admittedly, these a can be taken into account by adding
penalty terms to the loss/objective function, a technique often used
in the Physically Informed Neural Networks (PINNs) branch of research
(see \cite{CaiMaoWangYinKarn21,WangYuPerd22} for instance), but this
introduces new hyper-parameters needing calibration and does not
guarantee that constraints are strictly enforced. Moreover, a
sufficiently severe penalization of the constraints typically
deteriorates the problem's conditioning, possibly resulting in slow
convergence, especially if first-order methods are used.  This
suggests that there is a scope for a more direct approach.  The present
short paper proposes a step in that direction by considering adaptive
gradient algorithms for nonconvex problems which have bound 
constraints.

Nonconvex optimization with bound constraints is not a new subject,
but, to the author's knowledge, nearly all available algorithms use
objective function values to ensure global convergence by assessing
the effectiveness of a move from one iterate to the next (see, for
instance, \cite[Chapter~12]{ConnGoulToin00} for a description of suitable
trust-region methods and \cite[Section~14.1]{CartGoulToin22} for an
analysis of their complexity). This can be inconvenient, for instance in deep
learning applications where sampling techniques cause significant
noise (\cite{GratJeraToin23c} provides a convincing illustration of
this argument). Exceptions where the function values are not used in
the globalization strategy are the papers \cite{Chenetal19} and
\cite{AlacMaliCevh21}.  The first studies a stochastic zero-th order
(DFO) algorithm for problems with convex constraints based on the
Malahanobis metric where the gradient is approximated by
finite-differences and proves that, in expectation, the norm of the
projected gradient (in that metric) decreases like $k^{-1/4}$, where
$k$ is the iteration counter.  This study is motivated by applications
in adversarial training (see the references in \cite{Chenetal19} for
details).  The second paper considers a purely first-order stochastic
version of AMSgrad for weakly convex problems with convex constraints and
obtains a similar rate of global convergence. It also mentions
applications in adversarial training and reinforcement learning.
Although both approaches are based on projections, no specific mention
is made of the special case of bound constraints, despite the fact
that such constraints are an integral part of efficient formulations
of more complex problems with nonconvex constraints. Our aim in the
present paper is to complete this picture by considering a general
class of adaptive algorithms for general nonconvex problems including
the popular Adagrad  method, but also allowing the use of curvature
information, should it be available at a reasonable cost. The algorithm
we will introduce \emph{never evaluates the objective function} (as is
the case of the popular aforementioned Adagrad, ADAM and their
variants). We refer to this feature by the OFFO acronym, for Objective
Function Free Optimization. 

Our proposal is deterministic\footnote{It is the authors' opinion that
a study of the deterministic case is most useful for improving one's
understanding of an algorithm's behaviour in more general stochastic
contexts.} and is based on the combination of two existing algorithmic 
approaches: trust-region based projection
methods \cite[Chapter~12]{ConnGoulToin00} (more general and more
flexible than pure projected gradient methods) and a recent reinterpretation
of adaptive-gradient methods as a particular class of (the same)
trust-region techniques \cite{GratJeraToin22a}. We claim that this combination,
although only moderately complex technically, yields a useful class of
optimization algorithms, in particular given the interest of the
deep-learning community in methods capable of handling nonconvex
functions and bound constraints.

A related approach is that using variants of Frank-Wolfe algorithm
\cite{FranWolf56} to minimize a potentially nonconvex function subject
to convex constraints, whose complexity was analyed in
\cite{Laco16}. This method performs an objective-function-free
linesearch along the direction from the current iterate to the
minimizer of the first-order model in the feasible set. in this
context, the main selling point of the Frank-Wolfe algorithm has
always been (see \cite{FranWolf56}) that it replaces projection on the
feasible set (a quadratic optmization problem) by a potentially
cheaper linear programming subproblem whenever the feasible set is
defined by a set of linear inequalities. Further developments (see
\cite{Brauetal23}) have looked at stochastic versions of the method
replacing the linesearch by either a fixed stepsize rule of a
``function agnostic'' rule based on iteration number only. Some of
these variants (\cite{ZimmSpiePoku22} or \cite{Miaoetal22} and
references therein) make use of the definition of the feasible set as
a technique to bias the search direction towards (various forms of)
sparse solutions.  Our proposal differs from this line of work in that
we propose a ``function-aware'' adaptive stepsize strategy (the
trust-region mechanism, possibly using Adagrad-like adaptivity) using
first-order information more locally (and thus probably more
reliably), and that we also allow taking second-order information into
account when possible.  We also note that the (first-order)
trust-region based step also results from an extremely cheap linear
programming computation when the constraints are bounds on the
variables, thus avoiding the general quadratic optimization cost of a
projection in our case.

While adaptive-gradient methods for constrained nonconvex problems
appear to be little explored, this is not the case for problems featuring
a convex objective function, for which various methods have been proposed (see
\cite{BachLevy19,JoulRajGyorSzep20, EneNguyVlad21} for
instance). Unfortunately, these methods are difficult to adapt to the
nonconvex setting because they typically hinge on an acceleration
technique that, so far, require convexity.

Our short paper is organized as follows. We first propose a class of
OFFO adaptive-gradient methods for optimization with bound constraints
in Section~\ref{algo-s}. We then specialize this class in
Section~\ref{adag-s} to derive a suitably modified Adagrad algorithm
and analyze its evaluation complexity. We also propose, in
Section~\ref{divstep-s} another specialization of our class, in the
spirit of \cite{GratJeraToin22a}.   Some conclusions and perspectives are finally
discussed in Section~\ref{concl-s}.

\numsection{First-order minimization methods with bound constraints}
\label{algo-s}

The first optimization problem under consideration is given by
\beqn{problem}
\min_{x\in\calF} f(x)
\eeqn
where $f$ is a smooth function from $\Re^n$ to $\Re$ and
\beqn{calF-def}
\calF = \big\{ x \in \Re^n \mid \ell_i \leq x_i \leq u_i \tim{ for } i\in\ii{n} \big\}.
\eeqn
The component-wise lower bounds $\{\ell_i\}_{i=1}^n$ and upper bounds $\{u_i\}_{i=1}^n$
in \req{calF-def} satisfy $\ell_i \leq u_i$ for $i\in\ii{n}$. The
values $\ell_i= -\infty$ and $u_i=+\infty$ are allowed.
In what follows, we also assume the following:
\vspace*{2mm}
\noindent
\begin{description}
\item[AS.1:] the objective function $f(x)$ is continuously differentiable;
\item[AS.2:] its gradient $g(x) \eqdef \nabla_x^1f(x)$ is Lipschitz continuous with
   Lipschitz constant $L\geq 0$, that is
   \[
   \|g(x)-g(y)\| \le L \|x-y\|
   \]
   for all $x,y\in \Re^n$; 
\item[AS.3:] there exists a constant $\flow$ such that, for all $x\in\calF$, $f(x)\ge \flow$.
\end{description}

\noindent
AS.1, AS.2 and AS.3 are standard for the complexity analysis of optimization
methods seeking first-order critical points, AS.3 guaranteeing in particular
that the problem is well-posed. Note that we do not
assume the gradients to be uniformly bounded, at variance with
\cite{WuWardBott18,DefoBottBachUsun22,ZhouChenTangYangCaoGu20,GratJeraToin22a}. 

Because the purpose of the method we are about to describe is to compute approximate
first-order critical points for problem \req{problem}, we first review what we mean
by this statement. Given a vector $x\in\calF$, we first define
\beqn{chiixal-def}
\chi_i(x, \alpha) = \delta_i(x,\alpha) | g(x)_i |
\tim{ where }
\delta_i(x,\alpha)
= \max\big\{ \delta \in [0,\alpha] \mid x -\delta \,\sgn[g(x)_i]\, e_i \in \calF \big\},
\eeqn
with $g(x) = \nabla_x^1f(x)$ and $\alpha \geq 0$.
Observe that $\chi_i(x,\alpha)$ is the maximum
decrease of the linear model $g(x)^Ts$ achievable along the $i$-th coordinate vector
while preserving feasibility and the inequality $|s_i|\leq \alpha$, that is 
\beqn{chiixal-min}
\chi_i(x,\alpha)
= \big| \min\big\{ g(x)^Ts \tim{such that} s=-\delta \,\sgn(g_{i,k})\,e_i,\,
\delta \leq \alpha \tim{and} x+s \in \calF \big\}\big|
\eeqn
Hence $\chi_i(x,1) = |g_{i,k}|$ when $\calF = \Re^n$ (i.e., the problem
is unconstrained) or when the distance from $x$ to the bounds exceeds one
along its $i$-th coordinate. Moreover $\chi_i(x) \eqdef \chi_i(x,1)$ can be
interpreted as a continuous and backward stable first-order criticality
measure for the one-dimensional problem
\[
\min_{t \in \smallRe} \big\{ f(x+te_i) \tim{such that} \ell_i \leq x_i+te_i\leq u_i \big\}
\]
(see \cite{GratMoufToin11}). 
As it is standard to define an (unconstrained) $\epsilon$-approximate
first-order critical point at a point $x$ such that
\[
\|g(x)\| = \sqrt{\sum_{i=1}^n |g_i(x)|^2 } \leq \epsilon,
\]
we extend this notion here by defining a box-constrained
$\epsilon$-approximate first-order critical point as a vector $x\in
\calF$ such that 
\beqn{eps-critical}
\sqrt{\sum_{i=1}^n |\chi_i(x)|^2 } \leq \epsilon,
\eeqn
a backward stable condition whose left-hand side is
continuous.
If $\epsilon = 0$, \req{eps-critical}
defines an exact first-order critical point.

Of interest here are iterative methods which generate a sequence of
iterates $\{x_k\}_{k\geq0}$ where the step from $x_k$ to $x_{k+1}$ 
depends on the gradient of the objective function at $x_k$ and on algorithm-dependent \emph{weights}
$\{w_k = w(x_0, \ldots, x_k)\}$ whose main purpose is to control the step's size.
These weights are assumed to be bounded below by a strictly positive
constant, that is
\begin{description}
  \item[AS.4:] for each $i\in\ii{n}$ there exists a constant $\varsigma_i\in(0,1]$
    such that, $w_{i,k} \geq \varsigma_i$ for all $k\geq0$,
\end{description}

Given an iterate $x_k$, we add a subscript $k$ to and drop the 
argument $x_k$ from quantities of interest to indicate they are evaluated for
$x=x_k$, so that
\[
g_k = \nabla_x^1f(x_k), \ms
\delta_{i,k}(\alpha) = \delta_i(x_k,\alpha) \tim{ and }
\chi_{i,k} = \chi_i(x_k,1)
\]
for instance.

We now state the  \al{ASTR1B} algorithm (for Adaptively Scaled Trust Region using
1rst order information with Bounds), a variant of the \al{ASTR1}
method of \cite{GratJeraToin22a} capable of handling bound constraints.

\algo{ASTR1B}{\tal{ASTR1B}}
{
\begin{description}
\item[Step 0: Initialization. ]
  A starting point $x_0\in \calF$ is given. A constant $\tau \in
  (0,1]$ is also given.   Set $k=0$.
\item[Step 1: Define the trust-region. ]
  Compute $g_k = g(x_k)$ and define
  \beqn{Delta-def}
  \Delta_{i,k}  =\frac{\chi_{i,k}}{w_{i,k}}.
  \eeqn
\item[Step 2: Hessian approximation. ]
  Select a symmetric Hessian approximation $B_k$.
\item[Step 3: GCP.]
  Compute a step $s_k^L$ using
  \beqn{sL-def}
  s_{i,k}^L = - \sum_{i=1}^n \delta_{i,k}(\Delta_{i,k})\,\sgn(g_{i,k})
  \ms (i\in\ii{n})
  \eeqn
  and a step $s^Q_k$ given by
  \beqn{sQ-def}
  s^Q_k = \gamma_{k} s_k^L,
  \eeqn
  with
  \beqn{gamma-def}
  \gamma_k =
  \left\{ \begin{array}{ll}
  \min\left[ 1, \bigfrac{|g_k^Ts_k^L|}{(s_k^L)^T B_k s_k^L}\right] & \tim{if }
  (s_k^L)^T B_k s_k^L > 0,\\
  1 & \tim{otherwise.}
  \end{array}\right.
  \eeqn
  Finally select a step $s_k$ such that
  \beqn{sfeas}
  x_k + s_k \in \calF,
  \eeqn
  \beqn{sbound}
  |s_{i,k}| \le \Delta_{i,k}  \ms (i \in \ii{n}),
  \eeqn
  and
  \beqn{GCPcond}
  g_k^Ts_k + \half s_k^TB_ks_k \le \tau\left(g_k^Ts_k^Q + \half (s_k^Q)^TB_ks_k^Q\right).
  \eeqn

\item[Step 4: New iterate.]
  Define
  \beqn{xupdate}
  x_{k+1} = x_k + s_k,
  \eeqn
  increment $k$ by one and return to Step~1.
\end{description}
}

The \al{ASTR1B} algorithm, like \al{ASTR1}, allows the user to provide
(and use) second-order information in $B_k$, should it
be available. There is no need for $B_k$ to be the exact Hessian of
$f$ at $x_k$, and (limited-memory) quasi-Newton approximations are
acceptable. If $B_k$ is chosen to be indentically zero on every
iteration, the \al{ASTR1B} algorithm is then a ``purely first-order''
optimization method.

Albeit computing a step $s_k$ satisfying \req{sfeas}-\req{GCPcond} may seem
formidable at first sight, this is actually a simple task, and even more so if $B_k=0$. Indeed,
the step $s_k^L$ (which minimizes the linear model $g_k^Ts$ in the intersection
of the trust-region and the feasible set) is easily obtained from \req{sL-def}  since, using\req{chiixal-def},
\beqn{delta-2}
\delta_{i,k}(\Delta_{i,k}) = \left\{\begin{array}{ll}
\min\left[\Delta_{i,k},u_i - x_i\right] & \tim{if } g_{i,k} < 0\\
\min\left[\Delta_{i,k},x_i - \ell_i\right] & \tim{if } g_{i,k}> 0\\
0 & \tim{if } g_{i,k} = 0
\end{array}\right.
\eeqn
for $i\in\ii{n}$. If $B_k=0$, we
have that $s_k = s_k^Q = s_k^L$.  Otherwise the step  $s^Q_k$ is computed using \req{sQ-def}
and \req{gamma-def}, the latter possibly involving a single matrix vector product,
and may be seen a ``generalized Cauchy point'' (see \cite{GratJeraToin22a}), that
is a step minimizing the quadratic model along a good first-order
descent direction while preserving feasiblity with respect to the
bounds and the trust-region.  Any step preserving \req{sfeas}
and the bound \req{sbound}, and decreasing the value of the quadratic model
$g_k^Ts + \half s^TB_ks$ is then acceptable as $s_k$.  Such an improved step may be
computed for instance by using truncated projected Krylov techniques.
This ``generalized Cauchy point technique''  has been first proposed
in \cite{ConnGoulToin88a} and variants have been used extensively
since\footnote{In the non-OFFO context.}, for instance in the {\sf LANCELOT} \cite{ConnGoulToin92}, {\sf
  BOX-QUACAN} \cite{DiniGomeSant98b}, {\sf TRON} \cite{LinMore99} and
{\sf GALAHAD} \cite{GoulOrbaToin03b} packages.

As stated, the algorithm does not impose any restriction on $B_k$, but
it is clear that if curvature can be arbitrarily large, then the step
can be arbitrarily small and convergence may be impeded.  In the
convergence theory below, we therefore assume that
\begin{description}
  \item[AS.5:] There exists a constant $\kB\geq 1$
    such that $\|B_k\| \le \kB$ for all $k\geq 0$.
\end{description}

Also observe that, at variance with ``projected gradient methods''
using function values (see \cite[Section~12.2.1]{ConnGoulToin00}), the
\al{ASTR1B} algorithm does not use explicit projection onto the
feasible set, nor does it perform approximate
minimization of the quadratic model on the left-hand side of
\req{GCPcond} along the ``projected-gradient path'' when $B_k$ is nonzero. This has the
advantage, which we judge crucial in potential neural-network training
applications, of avoiding more than a single Hessian times vector
product,

Because there are many possible choices for the weights $w_{i,k}$ in the
\al{ASTR1B} algorithm, the latter may effectively be seen as a
\textit{class} of more specific methods, two of which we investigate
in Sections \ref{adag-s} and \ref{divstep-s}.

All results in this paper crucially depend on the following two
lemmas, which derive lower bounds on the decrease of the first-order
Taylor model of $f$ in the neighbourhood of an iterate $x_k$ and on
the value of $f$ itself. The first is a substantially modified version of
\cite[Lemma~12.2.2]{ConnGoulToin00}, while the second is an adaptation
of \cite[Lemma~2.1]{GratJeraToin22a}.

\llem{lindecratsL}{Suppose that AS.1 and AS.4 hold. Then, at each
  iteration $k\geq0$ generated by the  \al{ASTR1B} algorithm,
\beqn{ldsL}
|g_{i,k}s_{i,k}^L|\geq
\min\left[\frac{\chi_{i,k}^2}{w_{i,k}},\chi_{i,k}\right]
\ms ( i\in\ii{n})
\eeqn
and
\beqn{ldsL2}
|g_k^Ts_k^L|\geq \varsigma_{\min} \|s_k^L\|^2
\eeqn
where $\varsigma_{\min} = \min_{i\in\ii{n}}\varsigma_i$.
}

\proof{
 Let $i\in\ii{n}$ and observe that
\req{chiixal-def} and \req{sL-def} ensure that $s_{i,k}^L= 0$ and
$\chi_{i,k}=0$ whenever $g_{i,k}= 0$, so
that \req{ldsL} and \req{ldsL2} trivially hold in this case. We therefore
assume that $g_{i,k} \neq 0$ for the rest of the proof and
let $d_{i,k} = u_i-x_{i,k}$ if $g_{i,k} < 0$, or
$d_{i,k} = x_{i,k}-\ell_i$ if $g_{i,k}> 0$.  Note that \req{sL-def}
implies that
\beqn{sL-D}
|s_{i,k}^L| \leq \Delta_{i,k}
\eeqn
for all $i$ and $k$. \\

If $|s_{i,k}^L|\geq 1$, then $d_{i,k} \geq 1$ and $|g_{i,k}| = \chi_{i,k}$.  Moreover \req{sL-def}
implies that  $\Delta_{i,k} \geq 1$.  We thus obtain that
\beqn{dl-2}
|g_{i,k}s_{i,k}^L| \geq \chi_{i,k}\min[d_{i,k},\Delta_{i,k}] \geq \chi_{i,k}.
\eeqn
Suppose first that $|s_{i,k}^L| = \Delta_{i,k}$,  then, using
\req{Delta-def}, AS.4 and \req{sL-D},
\beqn{dl-1sg}
|g_{i,k}s_{i,k}^L|
= |g_{i,k}|\frac{\chi_{i,k}}{w_{i,k}}
=\frac{\chi_{i,k}^2}{w_{i,k}}
=w_{i,k}\frac{\chi_{i,k}^2}{w_{i,k}^2}
\geq  \varsigma_i\Delta_{i,k}^2
= \varsigma_i (s_{i,k}^L)^2.
\eeqn
If now $|s_{i,k}^L| = d_{i,k}$, one has $d_{i,k} \leq \frac{\chi_{i,k}}{w_{i,k}}$,  which yields
\beqn{dl-2sg}
|g_{i,k}s_{i,k}^L|
= |g_{i,k}|d_{i,k}
= \chi_{i,k}d_{i,k}
= w_{i,k}\frac{\chi_{i,k}}{w_{i,k}}d_{i,k}
\geq w_{i,k} d_{i,k}^2
\geq \varsigma_i (s_{i,k}^L)^2.
\eeqn

Suppose now that $|s_{i,k}^L|< 1$.  If $|s_{i,k}^L|= \Delta_{i,k}$
and $d_i \geq 1$, then again $|g_{i,k}| = \chi_{i,k}$
and, using \req{Delta-def},
\beqn{dl-3}
|g_{i,k}s_{i,k}^L|= \chi_{i,k}\Delta_{i,k} = \frac{\chi_{i,k}^2}{w_{i,k}},
\eeqn
while using \req{Delta-def}, AS.4 and \req{sL-D} now gives that
\beqn{dl-3b}
|g_{i,k}s_{i,k}^L| = \chi_{i,k}\Delta_{i,k}
= \frac{\chi_{i,k}^2}{w_{i,k}}
= w_{i,k}\frac{\chi_{i,k}^2}{w_{i,k}^2}
\geq \varsigma_i\Delta_{i,k}^2
\geq \varsigma_i (s_{i,k}^L)^2.
\eeqn
If $|s_{i,k}^L|= \Delta_{i,k}$ and $d_i< 1$, then $\chi_{i,k} =
|g_{i,k}| d_{i,k}$ and, using  \req{sL-D} again,
\beqn{dl-4}
|g_{i,k}s_{i,k}^L|
= |g_{i,k}|\Delta_{i,k}
= \chi_{i,k}\frac{\Delta_{i,k}}{d_i}
\ge \chi_{i,k}\Delta_{i,k}
=\frac{\chi_{i,k}^2}{w_{i,k}}.
\eeqn
We also deduce, using \req{Delta-def}, \req{sL-D} and AS.4, that, in
this case,
\beqn{dl-4b}
|g_{i,k}s_{i,k}^L|
= \frac{\chi_{i,k}^2}{w_{i,k}}
= w_{i,k}\frac{\chi_{i,k}^2}{w_{i,k}^2}
\geq \varsigma_i\Delta_{i,k}^2
\geq \varsigma_i (s_{i,k}^L)^2.
\eeqn

Finally, if $|s_{i,k}^L|< 1$ and $|s_{i,k}^L| < \Delta_{i,k}$, then
one of the bounds on variable $i$ must be active at $x_k+s_k^L$ and thus
\beqn{dl-5}
|g_{i,k}s_{i,k}^L| = \chi_{i,k},
\eeqn
and we obtain, using once more \req{Delta-def}, \req{sL-D}, AS.4 and
the inequality  $|s_{i,k}^L|< 1$, that
\beqn{dl-5b}
 |g_{i,k}s_{i,k}^L|
= w_{i,k}\frac{\chi_{i,k}}{w_{i,k}}
\geq \varsigma_i\Delta_{i,k}
\geq \varsigma_i |s_{i,k}^L|
\geq \varsigma_i (s_{i,k}^L)^2.
\eeqn
Combining \req{dl-2}, \req{dl-3}, \req{dl-4} and \req{dl-5} then yields \req{ldsL},
while combining \req{dl-1sg}, \req{dl-2sg}, \req{dl-3b}, \req{dl-4b}, \req{dl-5b} and
the inequality  gives that, for all $i$ and $k$,
\beqn{dl-6}
|g_{i,k}s_{i,k}^L|
\geq \varsigma_i (s_{i,k}^L)^2.
\eeqn
But \req{sL-def} implies that  $g_{i,k}s_{i,k}^L<0$ for all
$i$ and $k$, and thus that, for $k\geq 0$,
\[
|g_k^Ts_k^L|
=\left|-\sum_{i=1}^n |g_{i,k}s_{i,k}^L|\right|
= \sum_{i=1}^n |g_{i,k}s_{i,k}^L|.
\]
The inequality \req{ldsL2} then follows by summing \req{dl-6} for $i\in\ii{n}$ and
using the definition of $\varsigma_{\min}$.
}

\llem{lemma:GCP}{
Suppose that AS.1, AS.2, AS.4 and AS.5 hold. Then, at each iteration
$k\geq0$ generated by the \al{ASTR1B} algorithm,
\beqn{gen-decr}
f(x_{k+1})
\le f(x_k) - \frac{\tau\varsigma_{\min}}{2\kappa_B}
     \bigsum_{i=1}^n\min\left[\frac{\chi_{i,k}^2}{w_{i,k}}, \chi_{i,k}\right]
      + \half(\kB + L) \sum_{i=1}^n \frac{\chi_{i,k}^2}{w_{i,k}^2}.
\eeqn
}

\proof{
We now consider the quadratic model and suppose first that
$(s_k^L)^TB_ks_k^L>0$ and $\gamma_k < 1 $.  Then,
we deduce from \req{sQ-def}, \req{gamma-def}, \req{ldsL2} and AS.5 that
\beqn{DL-4}
\begin{array}{lcl}
g_k^Ts_k^Q+\half (s_k^Q)^TB_ks_k^Q
& = & \min_{\gamma}[g_k^T(\gamma s_k^L+\half (\gamma s_k^L)^TB_k(\gamma s_k^L)]\\*[2ex]
& = & -\bigfrac{(g_k^Ts_k^L)^2}{2(s_k^L)^TB_ks_k^L}\\*[2ex]
& \leq & - \bigfrac{\varsigma_{\min}}{2\kappa_B}|g_k^Ts_k^L|\\*[2ex]
\end{array}
\eeqn
If now 
$(s_j^L)^TB_js_j^L \le 0$ or $\gamma_j=1$, then \req{sQ-def}
and \req{sL-def} give that
\beqn{DL-4b}
g_k^Ts_k^Q+\half (s_k^Q)^TB_ks_k^Q
= g_k^Ts_k^L+\half (s_k^L)^TB_ks_k^L
\le \half g_k^Ts_k^L < 0
\eeqn
and \req{DL-4} then again follows from the bounds
$\kB \ge 1$ and $\varsigma_i\leq 1$ for $i\in\ii{n}$ (see AS.4).

Successively using AS.1--AS.2, \req{GCPcond}, \req{DL-4}
and \req{Delta-def} therefore yields that, for $j\ge0$,
\begin{align*}
f(x_{j+1})
& \le f(x_j) + g_j^Ts_j + \half s_j^TB_js_j - \half s_j^TB_js_j + \half L\|s_j\|^2 \label{posdef}\\*[1.5ex]
& \le f(x_j) + \tau\left(g_j^Ts_j^Q + \half (s_j^Q)^TB_js_j^Q\right)  + \half(\kB + L) \|s_j\|^2 \\
& \le f(x_j) - \frac{\tau\varsigma_{\min}}{2\kappa_B}
      \bigsum_{i=1}^n\min\left[\frac{\chi_{i,k}^2}{w_{i,k}}, \chi_{i,k}\right],
      + \half(\kB + L) \sum_{i=1}^n \Delta_{i,j}^2 \\
& \le f(x_j) - \frac{\tau\varsigma_{\min}}{2\kappa_B}
      \bigsum_{i=1}^n\min\left[\frac{\chi_{i,k}^2}{w_{i,k}}, \chi_{i,k}\right]
      + \half(\kB + L) \sum_{i=1}^n \frac{\chi_{i,k}^2}{w_{i,k}^2},
\end{align*}
which completes the proof.
} 

\numsection{Adagrad with bound constraints and second-order models}\label{adag-s}

In the unconstrained case, the well-known Adagrad algorithm
\cite{DuchHazaSing11} (in particular reframed as a trust-region method \cite{GratJeraToin22a})
uses weights given by $ w_{i,k} = \sqrt{\varsigma+\sum_{j=0}^k g_{i,k}^2}$.
We pursue the analogy mentioned above between $g_{i,k}$ in the unconstrained case and
$\chi_{i,k}$ in the box-constrained case by defining weights as follows.
For given
$\varsigma \in (0,1]$ and $\vartheta \in (0,1]$ let,
for all $i\in\ii{n}$ and for all $k \geq 0$,
\beqn{w-adag}
w_{i,k} \in \left[\sqrt{\vartheta}\, v_{i,k}, v_{i,k}\right]
\tim{where}
v_{i,k} \eqdef \left(\varsigma + \sum_{j=0}^k \chi_{i,j}^2\right)^\half.
\eeqn

The weights used by the well-know Adagrad algorithm are thus
recovered by setting $\vartheta=1$. \textit{When applied to unconstrained
problems, \al{ASTR1B} with \req{w-adag},
$\vartheta=1$ and $B_k=0$ is therefore identical to the (deterministic) Adagrad
method, and thus generalizes this method to the box-contrained case.}

Note that, for all $k\ge 0$ and all $i\in\ii{n}$, the first part of
\req{w-adag} implies that
\beqn{prop-adag}
\frac{\chi_{i,k}}{w_{i,k}} \leq \frac{1}{\sqrt{\vartheta}}
\tim{ and }
\min\left[\frac{\chi_{i,k}^2}{w_{i,k}}, \chi_{i,k}\right]
  \geq \sqrt{\vartheta}\frac{\chi_{i,k}^2}{w_{i,k}}
\eeqn
so that we may rewrite the bound \req{gen-decr} as
\beqn{gen-decr-2}
f(x_{k+1})
\le f(x_k) - \frac{\tau\vartheta\varsigma_{\min}}{2\kappa_B}
\bigsum_{i=1}^n\frac{\chi_{i,k}^2}{w_{i,k}}
+ \half(\kB + L) \sum_{i=1}^n \frac{\chi_{i,k}^2}{w_{i,k}^2},
\eeqn
from which we deduce, by summing over iterations 0 to $k$, that
\beqn{decr-adag}
 \frac{\tau\vartheta\varsigma}{2\kappa_B}
\sum_{j=0}^k\sum_{i=1}^n\frac{\chi_{i,j}^2}{w_{i,j}}
\leq f(x_0) - f(x_{k+1}) + \half(\kB + L)\sum_{j=0}^k \sum_{i=1}^n \frac{\chi_{i,j}^2}{w_{i,j}^2}.
\eeqn

\lthm{theorem:astr1b-adag}{Suppose that AS.1--AS.3 and AS.5 hold and that the
\al{ASTR1B} algorithm is applied to problem \req{problem} with its
'Adagrad-like' weights given by \req{w-adag}. Then
\beqn{gradbound}
\average_{j\in\iiz{k}}\sum_{i=1}^n \chi_{i,j}^2
\le \frac{\kap{adag}}{k+1},
\eeqn
with $\Gamma_0 \eqdef f(x_0)-\flow$ and
\beqn{kapadag-def}
\kap{adag} = \max\left\{
\varsigma,\bigfrac{1}{2}e^{\frac{2 \Gamma_0\vartheta}{n(\kappa_B + L)}},
\frac{1}{2\varsigma}
\left(\frac{8n\kB(\kB+L)}{\tau\vartheta^\sfrac{5}{2}}\right)^2
\,\left|W_{-1}\left(-\frac{\tau \varsigma\vartheta^\sfrac{5}{2}}{8n\kB(\kB+L)}\right)\right|^2
\right\},
\eeqn
where $W_{-1}$ is the second branch of the Lambert function \cite{Corletal96}.
}

\proof{
Given the inequality \req{decr-adag}, the proof of the theorem is a
variation on that of \cite[Theorem~3.2]{GratJeraToin22a}, where the $i$-th
component of the gradient $g_{i,k}$ is replaced by $\chi_{i,k}$,
the criticality measure for the $i$-th variable. It is detailed in the
Appendix for the sake of completeness.
} 

\noindent
As noted in \cite{GratJeraToin22a}, it is possible to give a weaker but more
 explicit bound on $\kap{adag}$ by using an upper
 bound on the value of the involved Lambert function. This can be obtained from
 \cite[Theorem~1]{Chat13} which states that, for $x>0$,
 \beqn{Lamb-bound}
 \left|W_{-1}(-e^{-x-1})\right| \leq 1+ \sqrt{2x} + x.
 \eeqn
 Remembering that, for $\gamma_1$ and $\gamma_2$ given by \req{ggu},
 $\log\left( \frac{\gamma_2}{\gamma_1}\right) \geq \log(3) > 1$
 and choosing $x = \log\left( \frac{\gamma_2}{\gamma_1}\right) -1 > 0$
 in \req{Lamb-bound} then yields that
 \beqn{Lamb-bound2}
 \left|W_{-1}\left(-\frac{\gamma_1}{\gamma_2} \right)\right|
 \le \log\left(\frac{\gamma_2}
     {\gamma_1}\right)+\sqrt{2\left(\log\left(\frac{\gamma_2}{\gamma_1}\right)-1\right)}.
 \eeqn
It is also possible to extend the definition of $s_k^L$ in \req{sL-def}
  by premultiplying it by a stepsize $\alpha_k\in[\alpha_{\min}, 1]$ for some
  $\alpha_{\min} \in (0,1]$. Our results again remain valid (with modified
constants).

Observe that, if the algorithm is terminated as soon as
\req{eps-critical} (as we argued in Section~\ref{algo-s}),  it must
stop at the latest at iteration 
  \beqn{eps-order}
  k = \kap{adag}^2  \epsilon^{-2}.
  \eeqn
This corresponds the the $\epsilon$-order $\calO(\epsilon^{-2})$, the
standard complexity order of first-order methods using function values
(see \cite[Chapter~2]{CartGoulToin22}).

It also results from \cite[Theorem~3.3]{GratJeraToin22a} (applied for
$\mu=\half$) that the complexity bound given by
Theorem~\ref{theorem:astr1b-adag} is essentially sharp (in the sense
of \cite{CartGoulToin18a}), because this is the case for its
unconstrained variant. More precisely, we have the following result.

\lthm{sharp1}{The bound \req{gradbound} 
 is essentially sharp in that, for each $\eta\in(0,1]$, there exists a univariate function $f_{\eta}$
  satisfying AS.1-AS.3 and AS.5 such that, when applied to minimize $f_{\eta}$
  without constraints from the
  origin, the \al{ASTR1B} algorithm with
  \req{w-adag}, $B_k=0$ and $\vartheta=1$ produce a sequence of 
  gradient norms given by $\chi_{1,0} = g_0 =-2$ and $\chi_{1,k}^2 =
  g_k^2 = \frac{1}{k^{1+2\eta}}$ for $k\geq1$.
  }

\proof{ See \cite[Theorem~3.3]{GratJeraToin22a} and note that
\[
\average_{j\in\iiz{k}}\chi_{1,k}^2
= \frac{4}{k+1} +\frac{1}{k+1}\sum_{j=1}^k\frac{1}{j^{1+2\eta}}
\leq \frac{4+\zeta(1+2\eta)}{k+1},
\]
where $\zeta(\cdot)$ is the Riemann zeta function, which is
well-defined and finite for arguments exceeding one. 
}
  
\numsection{A ``diminishing stepsizes'' variant}\label{divstep-s}

We assume, in this section, that
the weights $w_{i,k}$ are chosen such that, for some power
parameter $0< \nu \leq \mu < 1$, all $i\in\ii{n}$ and some constants
$\varsigma_i\in(0,1]$ and $\theta > 0$,
\beqn{w-divstep}
\max[\varsigma_i,v_{i,k}]\,  (k+1)^\nu \le w_{i,k} \le \max[\varsigma_i,v_{i,k}]\,(k+1)^\mu
\ms
(k \geq 0),
\eeqn
where, for each $i$, the $v_{i,k}$ satisfy the properties that
\beqn{vikprop}
v_{i,k+1} > v_{i,k}  \tim{ implies that } v_{i,k+1} \leq |\chi_{i,k+1}|
\eeqn
and
\beqn{viklow}
v_{i,k} \geq \sqrt{\vartheta}|\chi_{i,k}| 
\eeqn
for some $\vartheta \in (0,1]$.
The motivation for considering these alternative class of variants is
the interesting numerical performance \cite{GratJeraToin22a} of the
choice 
\[
v_{i,k} = \max_{j\in\iiz{k}}|\chi_{i,j}|
\]
which satisfies \req{vikprop} and \req{viklow}.
This choice of weights ensures that
\[
\frac{\chi_{i,j}}{w_{i,j}} \leq \frac{1}{\sqrt{\vartheta}}
\]
and, as in the previous section, \req{gen-decr} implies
\req{decr-adag}. We then obtain the following result.

\lthm{theorem:astr1b-divstep}{
Suppose that AS.1, AS.2, AS.3 and AS.5 hold and that the
\al{ASTR1B} algorithm is applied to problem \req{problem}, where
the weights $w_{i,k}$ are chosen in accordance with
\req{w-divstep}, \req{vikprop} and \req{viklow}.
Then, for any $\eta \in (0,\tau \vartheta\varsigma_{\min})$ and
\beqn{jstar-ming}
j_\eta \eqdef
\left(\frac{\kB(\kB+L)}{\varsigma_{\min}(\tau \varsigma_{\min}-\eta)}\right)^{\sfrac{1}{\nu}},
\eeqn
there exist a constant $\kappa_\diamond$, a subsequence
$\{k_\ell\}\subseteq \{k\}_{j_\eta+1}^\infty$ and an index
$k_\varsigma$ (where $\kappa_\diamond$ and $k_\varsigma$
only depend on the problem and the algorithmic constants) such that,
for all $k_\ell \ge k_\varsigma$,
\beqn{ngkbound-ming} 
\min_{j\in\iiz{k_\ell}}\sum_{i=1}^n \chi_{i,j}^2
\le \kappa_\diamond \frac{(k_\ell+1)^\mu}{k_\ell-j_\eta}
\le \frac{2\kappa_\diamond (j_\eta+1)}{k_\ell^{1-\mu}}.
\eeqn
}

\proof{
The proof is again a variation of the proof of \cite[Theorem~4.2]{GratJeraToin22a}
where the $i$-th component of the gradient $g_{i,k}$ is replaced by $\chi_{i,k}$,
the criticality measure for the $i$-th variable. The proof is also
simpler than that in \cite{GratJeraToin22a} because our choice of weights $w_{i,k}$
is slightly more restrictive. The details are once more given in the Appendix.
}

Because $j_\theta$ and $k_\varsigma$ only
depends on $\nu$ and problem's constants,
Theorem~\ref{theorem:astr1b-divstep} gives some indication on the
rate of convergence for iterations beyond an \emph{a priori}
computable iteration index. The formulation of the
theorem is nevertheless weaker than that of
Theorem~\ref{theorem:astr1b-adag} since \req{ngkbound-ming} only holds
for iterates along the subsequence $\{k_\ell\}$ and there is no
guarantee that the bound given by the right-hand-side is valid at
other iterations. Fortunately, the index $k_\ell$ in this right-hand
side is an index in the complete sequence of iterates, which does
not depend on the subsequence.
A stronger result not involving subsequences has been proved in the
unconstrained case under the stronger assumption that gradients
remain uniformly bounded \cite[Theorem~4.1]{GratJeraToin22b}, and
its extension to the bound constrained case is possible much in the
same fashion that Theorem~4.2 of \cite{GratJeraToin22a} has been adapted here.

When $\mu$ and $\nu$ tend to zero,  the $k$-order of convergence
beyond $j_\theta$ (as stated by \req{ngkbound-ming}) tends to $\calO(1/\sqrt{k_\ell})$,
which the order derived for the methods of the previous section and is the standard $k$-order
for first-order methods using evaluations of the objective function, albeit the value
of $j_\theta$ might increase.  Moreover this result is essentially sharp, as implied by the
following theorem.

\lthm{sharp2}{
The bound \req{ngkbound-ming} is essentially sharp in that, for any $\omega >
\half(1-\nu)$, there exists a univariate function $f_\omega(x)$
satisfying AS.1--AS.3 and AS.5
such that the \al{ASTR1B} algorithm with \req{w-divstep}, $\mu=\nu$
and $B_k=0$ applied to this function without constraints produces a sequence of
first-order criticality measures given by $\chi_{1,k} = \|g_k\| =
\frac{1}{(k+1)^\omega}$.
}

\proof{See \cite[Theorem~4.3]{GratJeraToin22a} and the comments in the proof of
  Theorem~\ref{sharp1} above.}

\numsection{Numerical illustration}

Without any ambition of completeness, we now illustrate the behaviour
of the \al{ASTR1B} algorithm on a small set of 22 bound-constrained
problems from the {\sf CUTEst} collection (as made available in Matlab 
through S2MPJ \cite{GratToin24}).  The problems under consideration are given
in Table~\ref{table:theproblems} in Appendix~\ref{setup-a}.  We compare six algorithms.
\begin{description}
\item[\al{ASTR1B(0)}:] the \al{ASTR1B} algorithm using \req{w-adag}
  and $B_k=0$ (i.e. momentumless Adagrad with bounds)
\item[\al{ASTR1B(1)}:] the \al{ASTR1B} algorithm using \req{w-adag}
  and $B_k$ given by a limited BFGS update \cite{LiuNoce89}
  using one secant pair;
\item[\al{ASTR1B(3)}:] the \al{ASTR1B} algorithm using \req{w-adag}
  and $B_k$ given by a limited BFGS update using three secant pairs;
\item[\al{TRInf(0)}:] the standard trust-region algorithm using an
  $\ell_\infty$ trust-region and a purely linear model ($B_k=0$);
\item[\al{TRInf(1)}:] the standard trust-region algorithm using an $\ell_\infty$ trust-region
 and $B_k$ given by a limited BFGS update using one secant pair;
\item[\al{TRInf(3)}:] the standard trust-region algorithm using an
  $\ell_\infty$ trust-region and $B_k$ given by a limited BFGS update using three secant pairs.
  
\end{description}
We also ran these algorithms  for an accuracy level
$\epsilon = 10^{-3}$ on our small illustrative problem set using 0\%, 5\%,
15\% and 25\% of relative Gaussian noise on the gradients (and
objective function for the trust-region algorithms), in order to explore their
sensitivity to random perturbations.  The details of our experimental
setup are presented in Appendix~\ref{setup-a}. Our choice of
considering the same level of noise on the objective-function value (when
used) and the gradient is motivated by ``finite-sum'' applications such as deep
learning where both the objective-function value and the gradient are
computed by sampling.  Although the authors are aware that better
results can be obtained for the trust-region algorithms if one is
ready to substantially increase the sample size for the objective function (see
\cite{BlanCartMeniSche19} or \cite{BellGuriMoriToin21b} for instance),
this is unnatural and considerably more expensive in the sampling context.

\begin{table} 
  \begin{center}
    \begin{tabular}{|l|c|c|c|c|c|c|c|}
      \hline
      &Noise     & \al{ASTR1B(0)}&\al{ASTR1B(1)}&\al{ASTR1B(3)}&\al{TRInf(0)}&\al{TRInf(1)}&\al{TRInf(3)}\\
      \hline
      Performance & 0\% & 0.54 & 0.48 &0.48 & 0.87 & 0.90 & 0.90 \\
      \hline
      Reliability & 0\% & 100.0\% & 95.5\%  & 95.5\%  & 90.9\% & 90.9\% & 90.9\% \\
                  & 1\% &  99.6\% & 100.0\% & 100.0\% & 5.5\%  &  5.9\% &  5.0\% \\
                  & 5\% & 100.0\% & 99.6\%  & 99.6\%  & 5.0 \% &  5.0\% &  5.5\% \\
                  & 15\%& 100.0\% & 100.0\% & 100.0\% & {\footnotesize not run} &  {\footnotesize not run} & {\footnotesize not run} \\
                  & 25\%&  99.1\% & 100.0\% & 98.2\%  & {\footnotesize not run} &  {\footnotesize not run} & {\footnotesize not run} \\
      \hline
    \end{tabular}
    \caption{\label{table:theresults}Performance and reliability of \al{ASTR1B}
      variants and corresponding trust-region algorithms as a function of relative noise}
  \end{center}
\end{table}

Our results are reported in Table~\ref{table:theresults}.  The
first line of the table gives the compared performance of the six
methods of interest in the absence of noise, measured as the area of
the method's curve in a performance profile counting the number of
iterations to convergence for the six methods and
for an abscissa between 0 and 10, divided by 10.  Hence the closest to
one, the best performance (see \cite{PorcToin21} or
\cite{GratJeraToin23c} for other uses of this synthetic measure).
Clearly, the trust-region methods using
function values outperform the \al{ASTR1B} variants, whose
performance could be qualified of mediocre. But the picture
changes completely when one considers the reliability of the
algorithms, reported in the remaining lines of the table as the
percentage of successfully solved problems.  From this point of
view, the reliability of the trust-region methods becomes extremely
poor as soon as some noise is added\footnote{A typical failure of
the trust-region approach occurs when noise of the objective-function
values results in an erratic ratio of achieved to predicted reduction,
itself causing a decrease of the trust-region radius to a point
where the algorithm is stalled.}, while that of the \al{ASTR1B}
remains remarkably constant.  This observation conforts a similar
conclusion obtained in \cite{GratJeraToin23c} for unconstrained problems.

\numsection{Conclusions}\label{concl-s}

We have combined existing optimization techniques to propose an
adaptive gradient algorithm for minimizing general nonconvex objective
functions subject to bound constraints and allowing the use of
curvature information.  We have also analyzed its
evaluation complexity and shown that it is (in order) identical to
that of both standard algorithms using function values for the same
problem and adaptive-gradient algorithms for the unconstrained one.
Interestingly, our theoretical results do not require weak-convexity.

We view this proposal as a useful step towards efficient first-order
OFFO algorithms for nonconvex problems with nonconvex constraints, a
subject of importance in ``constraint-aware'' machine learning
applications such as PINNs, adversarial training and other approaches.
Another interesting application of the ideas developped here is the
use of sparsity-inducing norms (such as those discussed in
\cite{Brauetal23}) to define the trust region geometry.

{\footnotesize

\section*{\footnotesize Acknowledgements}

}

\appendix
\setcounter{equation}{0}
\appnumsection{Proof of Theorem~\ref{theorem:astr1b-adag}}

The proof uses the following technical lemma, due to
\cite{WardWuBott19,DefoBottBachUsun22}.

\llem{gen:series}{Let $\{a_k\}_{k\ge 0}$ be a non-negative sequence,
$\xi>0$ and define, for each $k \geq 0$,
$b_k = \sum_{j=0}^k a_j$.  Then 
\beqn{alphasup1series-bound}
\sum_{j=0}^k  \frac{a_j}{(\xi+b_j)}
\le  \log\left(\frac{\xi + b_k}{\xi} \right).
\eeqn
}

\proof{See \cite{WardWuBott19},
  \cite{WuWardBott18} or 
  \cite[Lemma~3.1]{GratJeraToin22a}.
}

\textbf{Proof of  Theorem~\ref{theorem:astr1b-adag}}
We see from \req{w-adag} that $w_{i,k}$ verifies \textbf{AS.4} and we
may thus use Lemma~\ref{lemma:GCP} and its consequence \req{decr-adag}.

For each $i\in\ii{n}$, we apply Lemma~\ref{gen:series}
with $a_k = \chi_{i,k}^2$ and $\xi = \varsigma$ and obtain that,
\[
\bigsum_{i=1}^n\bigsum_{j=0}^k \frac{\chi_{i,k}^2}{w_{i,k}^2}
\le \frac{1}{\vartheta}\bigsum_{i=1}^n
\log\left(\frac{1}{\varsigma}\left(\varsigma+
\bigsum_{l=0}^k \chi_{i,l}^2\right)\right)
\le \frac{n}{\vartheta} \log\left(1 +  \frac{1}{\varsigma}\bigsum_{l=0}^k \sum_{i=1}^n \chi_{i,l}^2\right).
\]
and substituting this bound in \req{decr-adag} then gives that
\beqn{noname}
\frac{\tau\vartheta\varsigma}{2\kappa_B}
\sum_{j=0}^k\sum_{i=1}^n\frac{\chi_{i,j}^2}{w_{i,j}}
\leq \Gamma_0
+ \half(\kB + L)\frac{n}{\vartheta} \log\left(1 +  \frac{1}{\varsigma}\bigsum_{l=0}^k \sum_{i=1}^n \chi_{i,l}^2\right).
\eeqn
Suppose now that 
\beqn{logmuhalfhyp}
\sum_{j=0}^k \sum_{i=1}^n\chi_{i,j}^2
\geq \max\left[ \varsigma, \bigfrac{1}{2}e^{\frac{2 \vartheta \Gamma_0}{n(\kappa_B + L)}}\right],
\eeqn
\noindent
implying that 
\[
1 + \frac{1}{\varsigma}\sum_{j=0}^k  \sum_{i=1}^n\chi_{i,j}^2 \leq \frac{2}{\varsigma}\sum_{j=0}^k  \sum_{i=1}^n\chi_{i,j}^2
\tim{and}
\Gamma_0 \leq \frac{n(\kB + L)}{2\vartheta} \log\left( \frac{2}{\varsigma}\bigsum_{j=0}^k  \sum_{i=1}^n\chi_{i,j}^2\right).
\]
\noindent
Using \req{noname}, we obtain then that
\[
\bigfrac{\tau\sqrt{\varsigma}\vartheta^\sfrac{3}{2}}{2\sqrt{2}\,\kB \,\sqrt{\sum_{\ell=0}^k \sum_{i=1}^n\chi_{i,j}^2}} \bigsum_{j=0}^k  \sum_{i=1}^n\chi_{i,j}^2
\le \frac{n(\kB + L)}{\vartheta} \log\left(  \frac{2}{\varsigma} \bigsum_{j=0}^k \sum_{i=1}^n\chi_{i,j}^2 \right),
\]
that is
\beqn{tosolvmuhalf}
 \bigfrac{\tau \sqrt{2\varsigma} \vartheta^\sfrac{5}{2}}{4 \kB } \sqrt{  \bigsum_{j=0}^k \sum_{i=1}^n\chi_{i,j}^2}
\le 2 n(\kB + L) \log\left(  \sqrt{\frac{2}{\varsigma} \bigsum_{j=0}^k \sum_{i=1}^n\chi_{i,j}^2}
\right).
\eeqn
Now define
\beqn{ggu}
\gamma_1 \eqdef \bigfrac{\tau \varsigma\vartheta^\sfrac{5}{2} }{4 \kB },
\ms
\gamma_2 \eqdef  2 n(\kB + L)
\tim{ and }
u \eqdef \sqrt{\frac{2}{\varsigma} \bigsum_{j=0}^k\sum_{i=1}^n\chi_{i,j}^2}
\eeqn
and observe that that $\gamma_2 > 3 \gamma_1$ because $\tau\sqrt{\varsigma} \vartheta^\sfrac{3}{2}\leq 1$
and $\kB \geq 1$.
The inequality \req{tosolvmuhalf} can then be rewritten as
\beqn{tosolvemuhalf}
\gamma_1 u \le \gamma_2 \log(u).
\eeqn
Let us denote by $\psi(u) \eqdef \gamma_1 u - \gamma_2 \log(u)$. Since
$\gamma_2 > 3 \gamma_1$, the equation $\psi(u)=0$ admits two roots $u_1 \leq
u_2$ and \req{tosolvemuhalf} holds for $u\in[u_1,u_2]$.
The definition of $u_2$ then gives that
\[
\log(u_2)- \frac{\gamma_1}{\gamma_2}u_2 = 0
\]
which is
\[
u_2e^{-\frac{\gamma_1}{\gamma_2}u_2} = 1.
\]
Setting $z = -\frac{\gamma_1}{\gamma_2}u_2$, we obtain that
\[
z e^z = -\frac{\gamma_1}{\gamma_2}
\]
Thus $z = W_{-1}(-\frac{\gamma_1}{\gamma_2})<0$, where $W_{-1}$ is the second
branch of the Lambert function defined over $[-\frac{1}{e}, 0)$.
As $-\frac{\gamma_1}{\gamma_2} \geq -\frac{1}{3} $, $z$ is well defined and thus
\[
u_2
= -\frac{\gamma_2}{\gamma_1}\,z
= -\frac{\gamma_2}{\gamma_1}\,W_{-1}\left(-\frac{\gamma_1}{\gamma_2}\right)>0.
\]
As a consequence, we deduce from \req{tosolvemuhalf} and \req{ggu} that
\[
\bigsum_{j=0}^k\sum_{i=1}^n\chi_{i,j}^2
= \frac{\varsigma}{2}\,u_2^2
=\frac{1}{2\varsigma}\,\left(\frac{8n\kB(\kB+L)}{\tau\vartheta^\sfrac{5}{2}}\right)^2
\,\left|W_{-1}\left(-\frac{\tau \varsigma\vartheta^\sfrac{5}{2}}{8n \kB(\kB+L)}\right)\right|^2.
\]
and
\beqn{firstmuhalfineq}
\average_{j\in\iiz{k}} \sum_{i=1}^n\chi_{i,j}^2
\leq \frac{1}{2\varsigma}\,
\left(\frac{8n\kB(\kB+L)}{\tau\vartheta^\sfrac{5}{2}}\right)^2
\,\left|W_{-1}\left(-\frac{\tau \varsigma \vartheta^\sfrac{5}{2}}{8n
  \kB(\kB+L)}\right)\right|^2\cdot\frac{1}{k+1}.
\eeqn
If \eqref{logmuhalfhyp} does not hold, we have that
\beqn{secmuhalfineq}
\average_{j\in\iiz{k}} \sum_{i=1}^n\chi_{i,j}^2
< \max\left\{
\varsigma, \frac{1}{2}\, e^{\frac{2\Gamma_0\vartheta}{n(\kappa_B + L)}}\right\}\cdot \frac{1}{k+1}.
\eeqn
Combining \req{firstmuhalfineq} and \req{secmuhalfineq}
gives \req{gradbound}. \epr

\appnumsection{Proof of Theorem~\ref{theorem:astr1b-divstep}}

From \req{decr-adag} and AS.3, using $w_{\min,j}\eqdef \min_{i\in\ii{n}}w_{i,k}$ ensures that
\beqn{newb-ming}
\Gamma_0
\ge f(x_0)-f(x_{k+1})
\ge \sum_{j=0}^k \sum_{i=1}^n \frac{\chi_{i,j}^2}{2\kB  w_{i,j}}
    \left[\tau \vartheta\varsigma_{\min}-\frac{\kB(\kB+L)}{w_{\min,j} }\right].
\eeqn
Consider now an arbitrary $\eta \in (0,\tau \vartheta\varsigma_{\min})$ and suppose first
that, for some $j$, 
\beqn{brackle-ming}
\left[\tau \vartheta\varsigma_{\min}-\frac{\kappa_B+L}{w_{\min,j}}\right] \le \eta,
\eeqn
i.e., using \req{w-divstep},
\[
\varsigma_{\min}\,\,j^\nu \le w_{\min,j} \le \frac{\kB(\kB+L)}{\tau \vartheta\varsigma_{\min}-\eta}.
\]
But this is impossible for $j > j_\eta$ for $j_\eta$ given by \req{jstar-ming},
and hence \req{brackle-ming} fails for all $j > j_\eta$.
As a consequence, we have that, for $k> j_\eta$,
\begin{align}
f(x_{j_\eta+1}) - f(x_k)
&\ge \eta \sum_{j=j_\eta+1}^k \sum_{i=1}^n \frac{\chi_{i,j}^2}{2\kB  w_{i,j}}\neol
&\ge \frac{\eta}{2\kB} \sum_{j=j_\eta+1}^k \sum_{i=1}^n\frac{\chi_{i,j}^2}{\max[\varsigma_i,v_{i,j}]\, \theta\, (j+1)^\mu}\neol
&\ge \frac{\eta}{2\kB(k+1)^\mu }
     \sum_{j=j_\eta+1}^k \sum_{i=1}^n\min\left[\frac{\chi_{i,j}^2}{\varsigma_i},\frac{\chi_{i,j}^2}{v_{i,j}}\right]\neol
&\ge \frac{\eta (k-j_\eta)}{2\kB(k+1)^\mu}
     \min_{j\in\iibe{j_\eta+1}{k}}
     \left(\sum_{i=1}^n\min\left[\frac{\chi_{i,j}^2}{\varsigma_i},\frac{\chi_{i,j}^2}{v_{i,j}}\right]\right)\label{A1}
\end{align}
But we also know from \req{gen-decr}, \req{w-divstep} and \req{viklow} that
\begin{align}
f(x_0)-f(x_{j_\eta+1})
&\ge \sum_{j=0}^{j_\eta} \sum_{i=1}^n \frac{\tau\vartheta\varsigma_{\min}\chi_{i,j}^2}{2\kB w_{i,j}}
- \half\ (\kB+L)\sum_{j=0}^{j_\eta} \sum_{i=1}^n\frac{\chi_{i,j}^2}{w_{i,j}^2}\neol
&\ge - \half \kB(\kB+L)\sum_{j=0}^{j_\eta}\sum_{i=1}^n\frac{\chi_{i,j}^2}{w_{i,j}^2}\neol
&\ge - \frac{n\kB(\kB+L)}{2\vartheta}\,j_\eta\label{A2}
\end{align}
Combining \req{A1} and \req{A2}, we obtain that
\[
\Gamma_0
\ge f(x_0)-f(x_{k+1})
\ge - \frac{n\kB(\kB+L)}{2\vartheta}\,j_\eta + \frac{\eta (k-j_\eta)}{2\kB(k+1)^\mu}
     \min_{j\in\iibe{j_\eta+1}{k}}\left(\sum_{i=1}^n\min\left[\frac{\chi_{i,j}^2}{\varsigma_i},\frac{\chi_{i,j}^2}{v_{i,j}}\right]\right)
\]
and thus that
\[
\min_{j\in\iibe{j_\eta+1}{k}}\left(\sum_{i=1}^n\min\left[\frac{\chi_{i,j}^2}{\varsigma_i},\frac{\chi_{i,j}^2}{v_{i,j}}\right]\right)
\le \frac{2\kB(k+1)^\mu}{\eta (k-j_\eta)}\left[\Gamma_0 + \frac{n \kB(\kB+L)}{2\vartheta}\,j_\eta\right]
\]
and we deduce that there must exist a subsequence
$\{k_\ell\}\subseteq \{k\}_{j_\eta+1}^\infty$ such that, for each $\ell$,
\beqn{A3}
\sum_{i=1}^n\min\left[\frac{\chi_{i,k_\ell}^2}{\varsigma_i},\frac{\chi_{i,jk_\ell}^2}{v_{i,k_\ell}}\right]
\le \frac{2\kB(k_\ell+1)^\mu}{\eta (k_\ell-j_\eta)}\left[\Gamma_0 +\frac{n \kB(\kB+L)}{2\vartheta}\,j_\eta\right].
\eeqn
But
\beqn{A0}
\frac{(k_\ell+1)^\mu}{k_\ell-j_\eta}
< \frac{2^\mu k_\ell^\mu}{k_\ell-j_\eta}
< \frac{2 k_\ell^\mu}{k_\ell-j_\eta}
=\frac{2 k_\ell^\mu k_\ell}{(k_\ell-j_\eta)k_\ell}
=\frac{k_\ell}{k_\ell-j_\eta} \cdot \frac{2}{k_\ell^{1-\mu}}
\le  \frac{2(j_\eta+1)}{k_\ell^{1-\mu}},
\eeqn
where we used the facts that $\mu <1$ and that $\frac{k_\ell}{k_\ell-j_\theta} $
is a decreasing function for $k_\ell \ge j_\theta+1$.
Using this inequality, we thus obtain from \req{A3} that, for each $\ell$,
\[
\sum_{i=1}^n\min\left[\frac{\chi_{i,k_\ell}^2}{\varsigma_i},\frac{\chi_{i,k_\ell}^2}{v_{i,k_\ell}}\right]
\le \frac{4\kB(j_\eta+1)}{\eta \, k_\ell^{1-\mu}}\left[\Gamma_0 \theta
  + n \frac{\kB(\kB+L)}{2\vartheta}\,j_\eta\right].
\]
As a consequence,
\[
k_\varsigma \eqdef
\left(\frac{4\kB(j_\eta+1)\left[\Gamma_0 
      +\frac{n \kB(\kB+L)}{2\vartheta}\,j_\eta\right]}{\eta\varsigma_{\min}}\right)^\frac{1}{1-\mu}
\]
is such that, for all $k_\ell \ge k_\varsigma$,
\beqn{conseq}
\min\left[\frac{\chi_{i,k_\ell}^2}{\varsigma_i,},\frac{\chi_{i,k_\ell}^2}{v_{i,k\ell}}\right]
\le  \varsigma_{\min}.
\eeqn
But \req{viklow} ensures that
\[
\min\left[\frac{\chi_{i,k_\ell}^2}{\varsigma_i},\frac{\chi_{i,k_\ell}}{\sqrt{\vartheta}}\right]
 \leq \min\left[\frac{\chi_{i,k_\ell}^2}{\varsigma_i,},\frac{\chi_{i,k_\ell}^2}{v_{i,k\ell}}\right]
\le  \varsigma_{\min}.
\]
Now this inequality and the bounds
$\varsigma_i\in (0,1)$, $\vartheta \in (0,1)$ and $\varsigma_{\min}\in (0,1)$
together imply that $\chi_{i,k_\ell} \in (0,1)$ and hence that
$\chi_{i,k_\ell}^2 <\chi_{i,k_\ell}$. We thus obtain from \req{conseq}
that, for all $k_\ell \ge k_\varsigma$,
\begin{align*}
\sum_{i=1}^n \frac{\chi_{i,k_\ell}^2}{\max[\varsigma_i,\sqrt{\vartheta}]}
&\le \frac{2\kB(k_\ell+1)^\mu}{\eta (k_\ell-j_\eta)}\left[\Gamma_0 +\frac{n \kB(\kB+L)}{2\vartheta}\,j_\eta\right]
\end{align*}
which, because $\max[\varsigma_i,\sqrt{\vartheta}] \le 1$, gives that, for all $k_\ell \ge k_\varsigma$,
\beqn{choosetheta}
\sum_{i=1}^n\chi_{i,k_\ell}^2
\le
\frac{(k_\ell+1)^\mu}{k_\ell-j_\eta}\,
\left(\frac{4\kB}{\eta}\right)
  \left[\Gamma_0 + n \frac{\kB(\kB+L)}{\vartheta}\,j_\eta\right],
\eeqn
finally implying, because of  \req{A0}, that \req{ngkbound-ming} holds
with
\[
\kappa_\diamond = 2(j_\eta+1)\left(\frac{4\kB}{\eta}\right)
\left[\Gamma_0 + n \frac{\kB(\kB+L)}{\vartheta}j_\eta\right].
\]
\epr

\appnumsection{Details of the experimental setup}\label{setup-a}

The test problem used and their dimension are given in
Table~\ref{table:theproblems}.

\begin{table}[htb]
\begin{center}
\begin{tabular}{|lr|lr|lr|lr|}
\hline
BQPGABIM &       50 & HADAMALS &      400 & NCVXBQP1 &      500 & OBSTCLBL &      625 \\
BQPGASIM &       50 & JNLBRNG1 &      625 & NCVXBQP2 &      500 & OBSTCLBM &      625 \\
EXPLIN   &      600 & JNLBRNG2 &      625 & NCVXBQP3 &      500 & OBSTCLBU &      625 \\
EXPLIN2  &      600 & JNLBRNGA &      625 & NOBNDTOR &     1024 & QINGB    &      500 \\
EXPQUAD  &      120 & JNLBRNGB &      625 & OBSTCLAE &      625 && \\
GENROSEB &      500 & LINVERSE &      999 & OBSTCLAL &      625 && \\
\hline
\end{tabular}
\caption{\label{table:theproblems}The {\sf CUTEst}/S2MPJ problems used
for illustration}
\end{center}
\end{table}

Our test were performed in Matlab R2023b on a Dell Precision with 64GB of
memory and running Ubuntu 20.04. All algorithms were stopped after a
maximum number of 100000 iterations (remember we are mostly using
first-order methods).  All variants of the \al{ASTR1B} algorithms use
the constants given by
\[
\tau = 1, \ms
\vartheta = 1
\tim{ and }
\varsigma = 0.01.
\]
The trust-region algorithm was used with an acceptance threshold $\eta_1
= 10^{-4}$, a trust-region expansion threshold $\eta_2= 0.95$, a radius
contraction factor of $\half$ and a radius expansion
factor of $2$.  In all cases where a nonzero $B_k$ was used, the
approximate quadratic model was minimized using a truncated
Lanczos iteration with a maximum number of inner iterations given by
$3n$ and a relative gradient accuracy request of $10^{-4}$.
Whenever random noise is present in the evaluations, the gradient's and
(when relevant) the objective-function's values are perturbed according to
\[
f(x) = f(x)\big(1+\tau \calN(0,1)\big), \ms
[g(x)]_i = [g(x)]_i)\big(1+\tau \calN(0,1)\big), \ms (i\in\ii{n}),
\]
where $\tau \in [0,1]$ is the noise level and $\calN(0,1)$ is the
standard normal distribution. In these cases, the reported results are
based on the rounded averages of ten independent runs.

\end{document}